\documentclass[11pt]{article}
\usepackage{psfig,epsfig,latexsym,graphicx,here}

\input amssym.def
\input amssym.tex
\renewcommand{\thefootnote}{\fnsymbol{footnote}}

\def\slfrac#1#2{\hbox{\kern.1em %
 \raise.5ex\hbox{\the\scriptfont0 #1}\kern-.11em %
 /\kern-.15em\lower.25ex\hbox{\the\scriptfont0 #2}}}

\def\binom#1#2{{#1}\choose{#2}}

\newcommand{\sF}{{\cal F}}

\newcommand{\RR}{{\Bbb R}}

\newcommand{\dd}{\ldots}

\newcommand{\eeq}{\end{equation}}
\newcommand{\beql}[1]{\begin{equation}\label{#1}}
\newcommand{\eqn}[1]{(\ref{#1})}
\newcommand{\hsp}{\hspace*{\parindent}}

\makeatletter
\def\eqalignno#1{\displ@y \tabskip\@centering
  \halign to\displaywidth{\hfil$\@lign\displaystyle{##}$\tabskip\z@skip
    &$\@lign\displaystyle{{}##}$\hfil\tabskip\@centering
    &\llap{$\@lign##$}\tabskip\z@skip\crcr
    #1\crcr}}
\makeatother
\makeatletter
\def\@sect#1#2#3#4#5#6[#7]#8{\ifnum #2>\c@secnumdepth
     \def\@svsec{}\else
     \refstepcounter{#1}\edef\@svsec{\csname the#1\endcsname.\hskip .75em }\fi
     \@tempskipa #5\relax
      \ifdim \@tempskipa>\z@
        \begingroup #6\relax
          \@hangfrom{\hskip #3\relax\@svsec}{\interlinepenalty \@M #8\par}%
        \endgroup
       \csname #1mark\endcsname{#7}\addcontentsline
         {toc}{#1}{\ifnum #2>\c@secnumdepth \else
                      \protect\numberline{\csname the#1\endcsname}\fi
                    #7}\else
        \def\@svsechd{#6\hskip #3\@svsec #8\csname #1mark\endcsname
                      {#7}\addcontentsline
                           {toc}{#1}{\ifnum #2>\c@secnumdepth \else
                             \protect\numberline{\csname the#1\endcsname}\fi
                       #7}}\fi
     \@xsect{#5}}
\def\@begintheorem#1#2{\it \trivlist \item[\hskip \labelsep{\bf #1\ #2.}]}
\makeatother

\makeatletter
\def\tm{\raise5pt\hbox{{\rm\tiny TM}}}
\newbox\tmbox\setbox\tmbox=\hbox{\tm}%
\makeatother
\begin{document}
\begin{center}
{\Large {\bf McLaren's Improved Snub Cube and Other New Spherical Designs in Three Dimensions}} \\
\vspace{1\baselineskip}
{\em R. H. Hardin} and {\em N. J. A. Sloane} \\
\vspace{.25\baselineskip}
Mathematical Sciences Research Center \\
AT\&T Bell Laboratories \\
Murray Hill, NJ 07974 USA \\
\vspace{1.5\baselineskip}
September 11, 1995. Minor editorial changes July 23, 2002. \\
\vspace{2\baselineskip}
{\bf Abstract} \\
\vspace{.5\baselineskip}
\end{center}
\setlength{\baselineskip}{1.5\baselineskip}

Evidence is presented to suggest that, in three dimensions, spherical 6-designs with $N$ points exist for $N=24$, 26, $\ge 28$;
7-designs for $N=24$, 30, 32, 34, $\ge 36$;
8-designs for $N=36$, 40, 42, $\ge 44$;
9-designs for $N=48$, 50, 52, $\ge 54$;
10-designs for $N=60$, 62, $\ge 64$;
11-designs for $N=70$, 72, $\ge 74$;
and 12-designs for $N=84$, $\ge 86$.
The existence of some of these designs is established analytically, while
others are given by very accurate numerical coordinates.
The 24-point 7-design was first found by
McLaren in 1963, and --- although not identified as such by McLaren --- consists of the vertices of an ``improved'' snub cube,
obtained from Archimedes' regular snub cube (which is only a 3-design)
by slightly shrinking each square face and expanding each triangular
face.
5-designs with 23 and 25 points are presented which, taken together
with earlier work of Reznick, show that 5-designs exist for $N=12$, 16, 18, 20, $\ge 22$.
It is conjectured, albeit with decreasing confidence for $t \ge 9$, that these
lists of $t$-designs are complete and that no others exist.
One of the constructions gives a sequence of putative spherical $t$-designs
with $N= 12m$ points $(m \ge 2)$ where
$N= \frac{1}{2} t^2 (1+o(1))$ as $t \to \infty$.
\clearpage
\large\normalsize
\renewcommand{\baselinestretch}{1}
\thispagestyle{empty}
\setcounter{page}{1}
\begin{center}
{\Large {\bf McLaren's Improved Snub Cube and Other New Spherical Designs in Three Dimensions}} \\
\vspace{1\baselineskip}
{\em R. H. Hardin} and {\em N. J. A. Sloane} \\
\vspace{.25\baselineskip}
Mathematical Sciences Research Center \\
AT\&T Bell Laboratories \\
Murray Hill, NJ 07974 USA \\
\vspace{1.5\baselineskip}
\end{center}
\setlength{\baselineskip}{1.5\baselineskip}
\renewcommand{\thefootnote}{\arabic{footnote}}
\section{Introduction}
\hsp
A set of $N$ points $\wp = \{P_1 , \dd, P_N \}$ on the unit sphere $\Omega_d = S^{d-1} = \{ x= (x_1, \dd, x_d) \in \RR^d: x \cdot x =1 \}$ forms a
{\em spherical $t$-design} if the identity
\beql{eq1}
\int_{\Omega_d} f(x) d \mu (x) = \frac{1}{N} \sum_{i=1}^N
f(P_i)
\eeq
(where $\mu$ is uniform measure on $\Omega_d$ normalized
to have total measure 1) holds for all polynomials $f$ of degree $\le t$
(\cite{DGS}; \cite{GS81}; \cite[\S3.2]{SPLAG}).
In the present paper we are concerned only with the case $d=3$.\footnote{We are in the process of producing an analogous table of four-dimensional designs;
these will be described elsewhere.}

It is trivial that 1-designs exist if and only if $N \ge 2$, and Mimura \cite{Mim90} showed that 2-designs exist if and only if $N=4$,
$\ge 6$.
Bajnok \cite{Baj93} found 3-designs for $N=6$, 8, $\ge 10$ and conjectured
that they do not exist for $N=7$ and 9.\footnote{In his talk Bajnok actually claimed to have a 9-point 3-design, but he now believes that this was a mistake.}
That 3-designs do not exist for $N \le 5$ is a consequence of the lower bounds
$$
\eqalignno{
N & \ge \frac{(t+1)(t+3)}{4} ,~~{\rm if}~~ t~~{\rm odd} ~, & {\rm (2a)} \cr
N & \ge \frac{(t+2)^2}{4} ,~~{\rm if}~~ t ~~{\rm even} ~, & {\rm (2b)}\cr}
$$
$$
\begin{array}{l}
\mbox{if $t \neq 1,2,3,5$ the right-hand sides of (2a), (2b)} \\ [+.05in]
\mbox{can be increased by 1}~,
\end{array}
\eqno{{\rm (2c)}}
$$

\noindent
which were established in
\cite{DGS}, \cite{BD79}, \cite{BD80}.
In \cite{Me174} we showed that 4-designs exist for $N=12$, 14, $\ge 16$,
and conjectured that no others exist.
Reznick \cite{Rez95} showed that 5-designs exist for $N=12$, 16, 18, 20, 22, 24, $\ge 26$.
We have found 5-designs with $N=23$ and 25 (see Table~I), and, our
search having repeatedly failed in the remaining cases,
conjecture that 5-designs do not exist for $N=13$--15, 17, 19 and 21.
Bajnok \cite{Baj91} gave a general construction for $t$-designs on $\Omega_3$, but his designs (described in \S5) are much larger than ours.

Following Reznick \cite{Rez95}, we make use of the fact that a set of points $\{P_i\}$ forms a spherical $t$-design if and only if the polynomial identities
$$
\frac{1}{N} \sum_{i=1}^N (P_i \cdot x)^{2s} =
\left( \prod_{j=0}^{s-1} \frac{2j+1}{2j+3} \right) (x \cdot x)^s ~,
\eqno{{\rm (3a)}}
$$
and
$$
\frac{1}{N} \sum_{i=1}^N (P_i \cdot x)^{2 \overline{s} +1} = 0 ~,
\eqno{{\rm (3b)}}
$$
hold, where $s$ and $\overline{s}$ are defined by $\{ 2s, 2\overline{s} +1 \} = \{ t-1, t\}$ (see \cite{GS81}; \cite[p.~114]{Rez92}).
\section{Summary of results}
\hsp
Let $\tau (N)$ denote the largest value of $t$ for which an $N$-point
3-dimensional spherical $t$-design exists.
Since a $t$-design is also a $t'$-design for all $t' \le t$,
an $N$-point spherical $t$-design exists if and only if $\tau (N) \ge t$.

Our main results are summarized in Table~I, which gives what we believe are the values of $\tau (N)$ for $N \le 100$.
The assertions made in the first sentence of the abstract can then be
simply read off the table.
The table also gives, in columns 4 and 5, the largest symmetry group we have
found for such a design (using the notation of \cite{CM84}),
and in some cases a list of the sizes of the orbits under this group
and a description of the polyhedron formed by the points.
In most cases the designs found were not unique.

For every value of $N$ in the table we have found very accurate numerical coordinates
for a putative
spherical $t$-design with $t$ equal to the value given in column 2.
Furthermore, after a considerable amount of searching, we have been
unable to find a $(t+1)$-design, and so we conjecture that the entries in column 2 do indeed give the exact values of $\tau (N)$.

In a number of cases we have proved that there {\em is} a spherical $t$-design that is very close to our numerical approximation.
To do this we reduce (3a), (3b) to a set of simultaneous algebraic equations,
and then show either algebraically (with the help of Maple \cite{Map}) or
by interval methods (using Intbis \cite{Kear90}) that these
equations do have a solution in the neighborhood of the approximate solution.
Examples will be found in the next two sections.

A symbol V1 in the third column of Table~I indicates that we have an algebraic proof of the existence of the design,
V2 that we have a proof by interval methods, and V3 that we have a numerical solution with discrepancy $\Delta ( \wp )$ (defined below) at most $10^{-26}$.
References to the literature indicate who first proved the existence of some spherical $t$-design with this number of points (not necessarily the particular design described in the table).
\begin{table}[H]
\caption{Conjectured values of $\tau (N)$, the largest $t$ for which an $N$-point
configuration on the sphere in 3 dimensions forms a spherical $t$-design.}

\vspace*{+.1in}
{\small
\begin{tabular}{llllll}
\multicolumn{1}{c}{$N$} & \multicolumn{1}{c}{$\tau(N)$} & \multicolumn{1}{c}{Proof} & \multicolumn{1}{c}{Group} & \multicolumn{1}{c}{Order} & \multicolumn{1}{c}{Orbits (Description)} \\ [+.1in]
1 & 0 & V1 & $\infty$ & $\infty$ & 1 (single point) \\
2 & 1 & V1 & $\infty$ & $\infty$ & 2 (2 antipodal points) \\
3 & 1 & V1 & $[2,3]$ & 12 & 3 (equilateral triangle) \\
4 & 2 & V1 & $[3,3]$ & 24 & 4 (regular tetrahedron) \\
5 & 1 & V1 & $[2,3]$ & 12 & $3+2$ (triangular bipyramid) \\
6 & 3 & V1 & $[3,4]$ & 48 & 6 (regular octahedron) \\
7 & 2 & \cite{Mim90} & $[3]$ & 6 & $3^2+1$ \\
8 & 3 & V1 & $[3,4]$ & 48 & 8 (cube) \\
9 & 2 & \cite{Mim90} & $[2,3]$ & 12 & $6+3$ (triangular biprism) \\
10 & 3 & \cite{Baj93} & $[2^+,10]$ & 20 & 10 (pentagonal prism) \\
11 & 3 & \cite{Baj93} & $[2,3]^+$ & 6 & $6+3+2$ \\
12 & 5 & V1 & $[3,5]$ & 120 & 12 (regular icosahedron) \\
13 & 3 & \cite{Baj93} & $[4]$ & 8 & $4^3+1$ \\
14 & 4 & \cite{Me174} & $[2,3]^+$ & 6 & $6^2+2$ \\
15 & 3 & \cite{Baj93} & $[2,5]$ & 20 & $10+5$ \\
16 & 5 & \cite{Me174} & $[3,3]^+$ & 12 & $12+4$ (hexakis truncated tetrahedron) \\
17 & 4 & \cite{Me174} & $[2,3]^+$ & 6 & $6^2+3+2$ \\
18 & 5 & \cite{Rez95} & $[2^+,6]$ & 12 & $12+6$ \\
19 & 4 & \cite{Me174} & $[3]$ & 6 & $6^2+3^2+1$ \\
20 & 5 & V1 & $[3,5]$ & 120 & 20 (regular dodecahedron) \\
21 & 4 & \cite{Me174} & $[2,3]$ & 12 & $12+6+3$ \\
22 & 5 & \cite{Rez95} & $[2^+,10]$ & 20 & $10^2+2$ \\
23 & 5 & V2 & $[2,3]^+$ & 6 & $6^3+3+2$ \\
24 & 7 & \cite{McL63} & $[3,4]^+$ & 24 & 24 (improved snub cube) \\
25 & 5 & V1 & $[2,5]^+$ & 10 & $10^2+5$ \\
26 & 6 & V3 & $[2,3]^+$ & 6 & $6^4+2$ \\
27 & 5 & \cite{Rez95} & $[2,3]$ & 12 & $12^2+3$ \\
28 & 6 & V3 & $[2^+,4]$ & 8 & $8^3+4$ \\
29 & 6 & V3 & $[2]^+$ & 2 & $2^{14}+1$ \\
30 & 7 & V1 & $[3,4]^+$ & 24 & $24+6$ (tetrakis snub cube) \\
31 & 6 & V3 & $[5]^+$ & 5 & $5^6+1$ \\
32 & 7 & V1 & $[3,4]^+$ & 24 & $24+8$ (snub cube $+$ cube) \\
33 & 6 & V3 & $[2,3]^+$ & 6 &  \\
34 & 7 & V3 & $[2,4]^+$ & 8 &  \\
35 & 6 & V3 & $[2,5]^+$ & 10 & $10^3+5$ \\
36 & 8 & V3 & $[3,3]^+$ & 12 & $12^3$ (3 snub tetrahedra) \\
37 & 7 & V3 & $[3]^+$ & 3 &  \\
38 & 7 & V3 & $[3,4]^+$ & 24 & $24+8+6$ \\
39 & 7 & V3 & $[2,3]^+$ & 6 &  \\
40 & 8 & V3 & $[3,3]^+$ & 12 & $12^3+4$ \\
41 & 7 & V3 & $[2,3]^+$ & 6 &  \\
42 & 8 & V3 & $[2,4]^+$ & 8 &  \\
43 & 7 & V3 & $[6]^+$ & 6 &  \\
44 & 8 & V3 & $[3,3]^+$ & 12 & $12^3+4^2$ \\
45 & 8 & V3 & $[2]^+$ & 2 &  \\
46 & 8 & V3 & $[2,4]^+$ & 8 &  \\
47 & 8 & V3 & $[2,3]^+$ & 6 &  \\
48 & 9 & V1 & $[3,4]^+$ & 24 & $24^2$ (two snub cubes) \\
49 & 8 & V3 & $[4]^+$ & 4 &  \\
50 & 9 & V3 & $[2,6]^+$ & 12 & $12^4+2$
\end{tabular}
}
\end{table}
\setcounter{table}{0}
\begin{table}[H]
\caption{(cont.)~Conjectured values of $\tau (N)$, the largest $t$ for which an $N$-point
configuration on the sphere in 3 dimensions forms a spherical $t$-design.}

\vspace*{+.1in}
{\small
\begin{tabular}{llllll}
\multicolumn{1}{c}{$N$} & \multicolumn{1}{c}{$\tau(N)$} & \multicolumn{1}{c}{Proof} & \multicolumn{1}{c}{Group} & \multicolumn{1}{c}{Order} & \multicolumn{1}{c}{Orbits (Description)} \\ [+.1in]
51 & 8 & V3 & $[2,3]^+$ & 6 \\
52 & 9 & V3 & $[3,3]^+$ & 12 & $12^4+4$ \\
53 & 8 & V3 & $[2,3]^+$ & 6 \\
54 & 9 & V3 & $[3,4]^+$ & 24 & $24^2+6$ \\
55 & 9 & V3 & $[2]^+$ & 2 \\
56 & 9 & V3 & $[3^+,4]$ & 24 & $24^2+8$ \\
57 & 9 & V3 & $[2,3]^+$ & 6 \\
58 & 9 & V3 & $[2,4]^+$ & 8 \\
59 & 9 & V3 & $[2,3]^+$ & 6 \\
60 & 10 & V3 & $[3,3]^+$ & 12 & $12^5$ (5 snub tetrahedra) \\
61 & 9 & V3 & $[6]^+$ & 6 \\
62 & 10 & V3 & $[2,3]^+$ & 6 \\
63 & 9 & V3 & $[2,7]^+$ & 14  & $14^4+7$ \\
64 & 10 & V3 & $[3,3]^+$ & 12  & $12^5+4$ \\
65 & 10 & V3 & $[2]^+$ & 2 \\
66 & 10 & V3 & $[2,4]^+$ & 8 \\
67 & 10 & V3 & $[2]^+$ & 2 \\
68 & 10 & V3 & $[2^+,4]$ & 8 \\
69 & 10 & V3 & $[4]^+$ & 4 \\
70 & 11 & V3 & $[2,5]^+$ & 10  & $10^7$ \\
71 & 10 & V3 & $[2,3^+]$ & 6 \\
72 & 11 & V3 & $[3,5]^+$ & 60 & $60+12$ (pentakis truncated icosahedron) \\
73 & 10 & V3 & $[4]^+$ & 4 \\
74 & 11 & V3 & $[2,6]^+$ & 12 & $12^6+2$ \\
75 & 11 & V3 & $[2]^+$ & 2 \\
76 & 11 & V3 & $[3,3]^+$ & 12 & $12^6+4$ \\
77 & 11 & V3 & $[4]^+$ & 4 \\
78 & 11 & V3 & $[3,4]^+$ & 24  & $24^3+6$ \\
79 & 11 & V3 & $[2]^+$ & 2 \\
80 & 11 & V3 & $[3,5]^+$ & 60 & $60+20$ (hexakis truncated icosahedron) \\
81 & 11 & V3 & $[4]^+$ & 4 \\
82 & 11 & V3 & $[2^+,10^+]$ & 10 & $10^8+2$ \\
83 & 11 & V3 & $[2,3]^+$ & 6 &  \\
84 & 12 & V3 & $[3,3]^+$ & 12 & $12^7$ (7 snub tetrahedra) \\
85 & 11 & V3 & $[2,5]^+$ & 10 \\
86 & 12 & V3 & $[2,2]^+$ & 4 \\
87 & 12 & V3 & $[1]^+$ & 1 \\
88 & 12 & V3 & $[3,3]^+$ & 12 & $12^7+4$ \\
89 & 12 & V3 & $[2]^+$ & 2 \\
90 & 12 & V3 & $[2,4]^+$ & 8 \\
91 & 12 & V3 & $[2]^+$ & 2 \\
92 & 12 & V3 & $[3,3]^+$ & 12 & $12^7+4^2$ \\
93 & 12 & V3 & $[4]^+$ & 4 \\
94 & 13 & V3 & $[2^+,2^+]$ & 2 \\
95 & 12 & V3 & $[2]^+$ & 2 \\
96 & 13 & V3 & $[3,3]^+$ & 12 & $12^8$ (8 snub tetrahedra) \\
97 & 12 & V3 & $[4]^+$ & 4 \\
98 & 13 & V3 & $[2,4]^+$ & 8 \\
99 & 12 & V3 & $[2]$ & 4 \\
100 & 13 & V3 & $[3,3]^+$ & 12 & $12^8+4$
\end{tabular}
}
\end{table}
\clearpage

The numerical coordinates for these $t$-designs were found by a modified version of the
Hooke and Jeeves \cite{HJ61}
``pattern search'' optimizer that we have already used to search for spherical codes \cite{HSS}
and experimental designs \cite{Goss},
\cite{Me174}, \cite{JSPI}.
Let $\sF_d~(0 \le d \le t )$ denote the set of
${\binom{d+2}{2}}$
monomials $f = x_1^{e_1} x_2^{e_2} x_3^{e_3}$ of degree $d$,
and let $\Delta_f ( \wp )$ be the difference between the right and
left sides of \eqn{eq1} for this $f$ for a set of points
$\wp = \{ P_1, \ldots, P_N \}$.
The criterion we used was to minimize
$$\sum_{f \in \sF_{t-1}} \frac{(t-1)!}{e_1 ! e_2 ! e_3 !} \Delta_f (\wp)^2 +
\sum_{f \in \sF_t} \frac{t!}{e_1 ! e_2 ! e_3 !}
\Delta_f ( \wp )^2 ~,
$$
since $\wp$ is a spherical $t$-design if and only if this sum vanishes.
(The multinomial coefficients make the sums rotationally invariant.)
As a check we also computed the {\em discrepancy} of the points,
$$\Delta ( \wp ) = \sum_{d=1}^t \sum_{f \in \sF_t} \Delta_f ( \wp )^2 ~.$$
In practice we have found that in the range of Table~I, if $\Delta ( \wp ) < 10^{-16}$ then there is a spherical $t$-design $\wp^\ast$ very close to $\wp$.

Incidentally, it would be nice to have a formal version of this:
a theorem giving an explicit bound $\Delta_0 (t,N)$ such that if
a set of $N$ points in $\Omega_3$ satisfies $\Delta ( \wp ) < \Delta_0 (t,N)$, then a spherical $t$-design $\wp^\ast$ exists near $\wp$.

The search was conducted by choosing a symmetry group from the lists of
decomposable rotation groups of orders up to 21 and all indecomposable rotation
groups,
picking a random starting configuration invariant under this group,
and optimizing with respect to the above criterion in such a way as to
preserve (or increase) the symmetry.
The program cycled through the values of $N$ from 10 to 100.
For each $N$, equal effort was spent in trying to increase the value of $t$,
and in trying to find a larger group for the current $t$.
The search was terminated when no further improvements
were found after several months of computing.	

In Sections~3 and 4 we describe in more detail several of the designs mentioned in Table~I.
Numerical coordinates for all these designs
have been placed on NJAS's home page.
\section{The improved snub cube}
\hsp
The regular snub cube (\cite{Cox20}, \cite{Cun1}), the familiar
Archimedean solid with equal edges, has symmetry group
$[3,4]^+$.
We take this group to consist of all even permutations of the three
coordinates combined with any even number of sign changes,
and all odd permutations combined with any odd number of sign changes.
Then the vertices of the regular snub cube consist of the point $(A,B,C)$ and
its images under the group, where $A= .8503 \dd$, $B= .4623 \dd$,
$C= .2514 \dd$ are the unique positive roots of
$$\begin{array}{l}
7Z^6 + Z^4 - 3Z^2 -1 ~, \\ [+.05in]
7Z^6 - 3Z^4 + 5Z^2 -1 ~, \\ [+.05in]
7Z^6 - 19 Z^4 + 17 Z^2 -1 ~,
\end{array}
$$
respectively, and $A^2 + B^2 + C^2 =1$.
One may verify from Eq.~(2) that these 24 points form a spherical 3-design but
not a 4-design.

However, by moving the vertices slightly, we can obtain a 7-design.
Again we take the vertices to consist of the 24 images of $(A,B,C)$ under the
group, where $A^2 + B^2 + C^2 =1$.
Eq.~(3b) with $\overline{s} = 3$ is trivially satisfied, and Eq.~(3a) with $s=3$ leads to the equations $A^6 + B^6 +C^6 = 3/7$,
$A^4 (B^2 +C^2) + B^4 (C^2 + A^2) + C^4 (A^2 + B^2) = 6/35$,
$A^2 B^2 C^2 = 1/105$.
It is easy to show that
these equations are satisfied by taking $A= .86624682 \dd$,
$B= .42251865 \dd$, $C= .26663540 \dd$ to be the positive
roots of the single equation
$$105 Z^6 - 105 Z^4 + 21 Z^2 -1~.$$
The convex hull of these points is the ``improved'' snub cube, differing from the regular one in that each square face has been slightly shrunk and each triangular face slightly expanded.
It is almost indistinguishable in appearance from the regular snub cube.
As far as we know this polyhedron is new.\footnote{(Added later.)
The polyhedron still seems to be new.
However, Bruce Reznick has pointed out to us that this spherical 7-design was first
found by McLaren in 1963 (\cite{McL63}; \cite{Str71}, p.~298;
\cite{Rez92}, pp.~112--113).
This is a very nice design, and we give it up grudgingly.}

In 1981 Goethals and Seidel \cite{GS81a} had shown that a similar
improvement can be made to the regular
truncated icosahedron (or soccer ball), another of the Archimedean solids.
The 60 vertices of the Archimedean solid form a spherical 5-design, but
Goethals and Seidel showed that a slight perturbation of the vertices
(while preserving the group) changes them to a 9-design.
Again the improved version is almost indistinguishable from the original,
which is shown in Fig.~1(a).
However, as can be seen from Table~I, it is possible to find
a 9-design with only 48 points, and a 10-design with 60 points.
The convex hull of our 60-point 10-design is shown in Fig.~1(b).
Coordinates will be found in \S6.
This figure has symmetry group $[3,3]^+$, and is the union of five
snub tetrahedra.
It has 174 edges and 116 triangular faces, and we do not expect it to replace the standard soccer ball!
\begin{figure}[htb]
\centerline{(a)\includegraphics[angle=270,
width=2.5in]{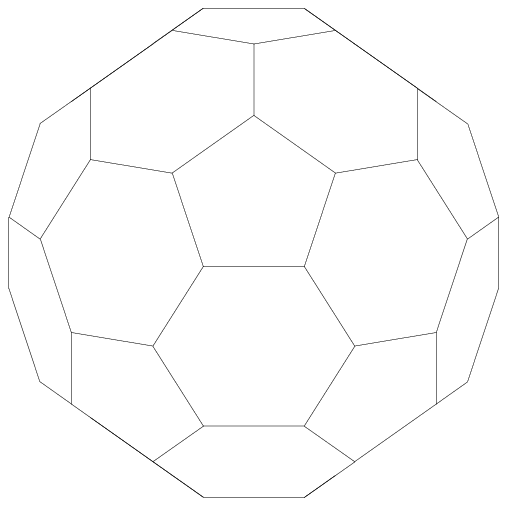}}
\centerline{(b)\includegraphics[angle=270,
width=2.5in]{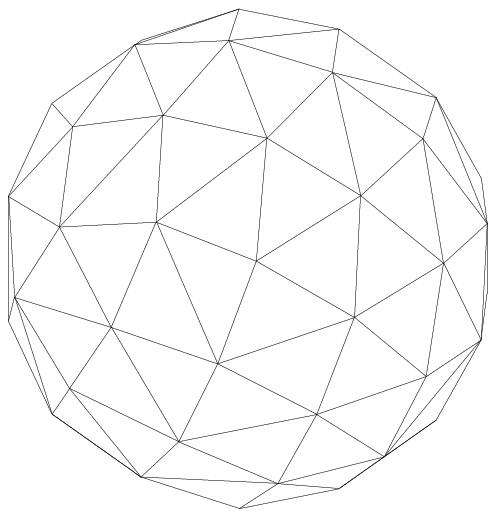}}
\caption{(a)~The regular truncated icosahedron (or soccer ball),
whose 60 vertices form a spherical 5-design.
The Goethals-Seidel improved football \protect\cite{GS81a},
which forms a spherical 9-design, is almost indistinguishable
from this.
(b)~Our 60-point spherical 10-design.}
\end{figure}
\section{Other examples of new spherical designs}
\hsp
We begin with two 5-designs that Reznick \cite{Rez95} was not able to find.
As one might expect, these are somewhat complicated.

{\em A 25-point 5-design with group $[2,5]^+$ of order 10.}
There are infinitely many 25-point 5-designs, of which the following
is the nicest we have found.
The points are
$$
\begin{array}{ccc}
(0 & \cos\,k \theta & \sin \, k \theta ) ~, \\ [+.1in]
(h_1 & - g_1 \,\cos\,k \theta & - g_1 \,\sin\,k \theta ) ~, \\ [+.1in]
(-h_1 & - g_1 \,\cos\,k \theta & g_1 \,\sin\,k \theta ) ~, \\ [+.1in]
(h_2 & g_2 \,\cos (k \theta + \pi_2 ) & g_2 \, \sin (k \theta + \pi_2 )) ~, \\ [+.1in]
(-h_2 & g_2 \,\cos (k \theta + \pi_2 ) & -g_2 \,\sin (k \theta + \pi_2)) ~,
\end{array}
\eqno{{\rm (4)}}
$$
where $0 \le k \le 4$, $\theta = 2 \pi /5$,
$$g_1 , g_2 = \frac{1}{2 \sqrt{3}}
\sqrt{7 \mp \sqrt{11}} = .5540 \dd , ~
.9272 \dd ~,$$
$h_1 = \sqrt{1-g_1^2} = .8325 \dd$,
$h_2 = \sqrt{1-g_2^2} = .3745 \dd$,
$\pi_2 = 2.057 \dd$ radians, defined by the condition that $\cos (\pi_2) = - .4670 \dd$ is a root of
$$16Z^5 - 20 Z^3 + 5Z -
\left( \frac{30 - 7 \sqrt{11}}{19} \right)^{5/2} +
\frac{1}{2}
\left( \frac{6 (7- \sqrt{11})}{19} \right)^{5/2} ~.
$$
(It is straightforward to show that these values satisfy
the equations obtained when (4) is substituted in (3a) and (3b) with $s= \overline{s} =2$.)
Other solutions can be obtained by including a phase
angle $\pi_1$ in the second and third lines of (4).

{\em A 23-point 5-design with a group of order 2.}
We must satisfy Eqs. (3a), (3b) with $s= \overline{s} =2$.
After a considerable amount of experimenting we found a
numerical solution with a symmetry of order 2, consisting of the points
\setcounter{equation}{4}
\beql{eq5}
(\pm 1,0,0) , (0, \pm 1, 0), (0, 0, -1),
(\pm a_i , \pm b_i , c_i) ~,
\eeq
where $0 \le i \le 8$, the $\pm 1$ signs in the last expression are linked,
$a_0 = 1/3$,
\beql{eq6}
a_i^2 + b_i^2 + c_i^2 = 1 , ~~~~0 \le i \le 8 ~,
\eeq
and the approximate values of the 26 unknowns
$a_1, \dd,$ $a_8 , b_0 , \dd,$ $b_8 , c_0 , \dd, c_8$ are
\beql{eq7}
\begin{array}{l}
.5654, .1949, .8837, .6521, .5610, .7414, .1927, .5854, -.2194, \\ [+.05in]
.3485, -.7812, -.4754, .7082, -.7301, .4805, .7199, -.2092, -.9169 \\ [+.05in]
 -.7476, -.5931, -.2807, -.2705, .3903, .4685, .6668, .7833
\end{array}
\eeq
respectively (only 4 decimal places are shown, although we worked with 12 places).
Substituting the symbolic forms \eqn{eq5} (with $a_0 = 1/3$) into (3a), (3b) with $s= \overline{s} =2$ produces 17 further equations, a typical one being
$$\frac{1}{9} c_0^2 + \sum_{i=1}^8 a_i^2 c_i^2 = \frac{23}{30} ~.$$
We then used interval Newton methods, as implemented in the software
package Intbis \cite{Kear90}, to show that these 26 equations have a unique solution in a small box around the point \eqn{eq7}.
We later found numerical solution with a larger group, $[2,3]^+$, of order 6 (see Table~I), but we have included the above existence proof as illustrative of the interval method.

{\em 30- and 32-point 7-designs with group $[3,4]^+$.}
These are similar to the improved snub cube described in \S3.
For 30 points we take the $24+6$ images of the points $(A,B,C)$, $(1,0,0)$ under the group, where $A= .7980 \dd$, $B= .5488 \dd$, $C= .2492 \dd$ are the positive roots of $84 Z^6 - 84Z^4 + 21 Z^2 -1$.
For 32 points we take the $24+8$ images of $(A,B,C)$, $\left( \frac{1}{\sqrt{3}}, \frac{1}{\sqrt{3}} , \frac{1}{\sqrt{3}} \right)$, where
$A= .8989 \dd$, $B= .4355 \dd$, $C= .0480 \dd$ are the positive roots of
$2835 Z^6 - 2835 Z^4 + 441 Z^2 -1$.
In both cases it is easy to show using Maple \cite{Map}
(in particular its Gr\"{o}bner basis package) that
equations (3a) and (3b) are satisfied.

{\em A 48-point 9-design with group $[3,4]^+$.}
Similar to the previous examples, but now we take the images of $(A,B,C)$, $(D,E,F)$, where $A= .9334 \dd$, $B= .3535 \dd$, $C= -.0620 \dd$, $D= .7068 \dd$,
$E= .6397 \dd$, $F= .3018 \dd$ are roots of the irreducible polynomial
$$
\begin{array}{l}
8141081016796875\,Z^{36}-48846486100781250\,Z^{34}+131885512472109375 \,Z^{32} \\ [+.1in]
~~~~~~-212133311066250000\,Z^{30}+226833777359437500\,Z^{28}-170368273215000000\,Z^{26} \\ [+.1in]
~~~~~~+92508869648475000\,Z^{24}-36735117403950000 \,Z^{22}+10602550092251250\,Z^{20} \\ [+.1in]
~~~~~~-2145915231232500\,Z^{18}+ 270106833039750\,Z^{16}-9766335726000\,Z^{14} \\ [+.1in]
~~~~~~-3473862884100\,Z^{12}+ 770657554800\,Z^{10}-80424958320\,Z^{8} \\ [+.1in]
~~~~~~+4880358000\,Z^{6}-168429429\,Z ^{4}+2733318\,Z^{2}-8269
\end{array}
$$
The complexity of this polynomial indicates why we have been satisfied to find purely numerical solutions for the larger designs in the table.
\section{Designs with larger numbers of points}
\hsp
Although Table~I only extends to $N=100$, larger designs for fixed $t$
may be obtained using the fact that an $N_1$-point design and an $N_2$-point
design can be combined to form an $(N_1 + N_2)$-point design.
For example $N$-point 6-designs can be found for all $N \ge 28$ by combining
the designs in the table.

Alternative (and exact) designs can be found using a construction of
Bajnok \cite{Baj91}.
An $n$-point {\em interval $t$-design} consists of $n$
distinct points $P_1 , \dd , P_n$ with $-1 \le P_i \le 1$ such that
$$\frac{1}{2} \int_{-1}^1 f(x) dx =
\frac{1}{n} \sum_{i=1}^n f(P_i)$$
holds for all polynomials $f$ of degree $\le t$.
Bajnok shows that by taking regular $m$-gons at latitudes
$P_1, \dd, P_n$ one obtains a 3-dimensional $mn$-point spherical $t$-design,
provided $m \ge t +1$.

It is known (see the survey by Gautschi \cite{Gau75}) that for $t=1,2,3, \dd , 11$, $n$-point interval $t$-designs exist for all $n \ge 1$, 2, 2, 4, 4, 6, 6, 9, 9, 13, 13, respectively.
When $t=6$, for example, Bajnok's construction produces $N$-point spherical 6-designs with $N=42$, 48, 49, 54, $\dd$ and all $N \ge 108$.
 
\section{A conjecturally infinite family of $t$-designs}
\hsp
Inspection of Table~I shows that there is a sequence of $N= 12m$-point spherical $t$-designs with group $G= [3,3]^+$ (or larger) which for $m=2, \dd, 8$ have $t=7, \dd, 13$.
One might naively expect this sequence to continue in a linear fashion,
but the true situation is more complicated.

A full orbit under $G$ can be taken to consist of 12 points
$(\pm A, \pm B, \pm C)$, $(\pm B , \pm C, \pm A)$, $(\pm C, \pm A, \pm B)$,
where the product of the signs is positive and $A^2 + B^2 + C^2 =1$.
(Their convex hull is a snub tetrahedron.)
So a set $\wp$ which is the union of $m$ full orbits under $G$ contains $2m$ degrees
of freedom.

Consider the ring $R$ of polynomials in $X$, $Y$, $Z$ that are
invariant under $G$, ignoring the trivial invariant $X^2 + Y^2 +Z^2$.
If $R_j$ is the subspace of $R$ consisting of homogeneous invariants of
degree $j$,
then the dimensions $d_j = \dim R_j$ are given by
the Molien series for $G$:
\beql{eq8}
\Phi ( \lambda ) = \frac{1+ \lambda^6}{(1- \lambda^3) ( 1- \lambda^4)}
= \sum_{j=0}^\infty d_j \lambda^j
\eeq
(see \cite{CM84}, Table~10).

In order for $\wp$ to form a $t$-design it is necessary and sufficient that the average
of $f$ over $\wp$ is equal to the average of $f$ over $\Omega_3$ for all
$f \in R_1 \cup R_2 \cup \cdots \cup R_t$ (\cite{GS79,GS81a}).
This imposes
\beql{eq9}
e_t = d_1 + d_2 + \cdots d_t
\eeq
conditions on $\wp$.
So provided $2m \ge e_t$, we may reasonably expect that it will be possible to choose the orbits
so that all the conditions are satisfied, and then a $t$-design with $N=12m$ points will exist.
The values of $e_t$ can be obtained by expanding \eqn{eq8}, and we discover that $t$-designs with $N=12m$ points should exist for the values of $t$ and $N$ shown in Table~II.
Table~I shows that such designs do indeed exist for $t \le 13$
(in fact for $t=7$ only 24 points are needed).
We have verified numerically that the predicted designs also exist for all $t \le 21$, and Table~III gives a set of orbit representatives
$(A,B,C)$ for a selection of these designs.
(The others can be obtained from NJAS's home page --- see \S2.)
\begin{table}[H]
\caption{Beginning of conjecturally infinite sequence of 3-dimensional spherical $t$-designs with $N=12m$ points and group $[3,3]^+$.
These have been constructed numerically for $t \le 21$ (cf. Table~III).}
$$
\begin{array}{rrrrrrrrrrrr}
t: & 3 & 4 & 5 & 6 & 7 & 8 & 9 & 10 & 11 & 12 & 13 \\
N: & 12 & 12 & 12 & 24 & 36 & 36 & 48 & 60 & 72 & 84 & 96 \\
~ \\
t: & 14 & 15 & 16 & 17 & 18 & 19 & 20 & 21 & 22 & 23 & 24 \\
N: & 108 & 132 & 144 & 156 & 180 & 204 & 216 & 240 & 264 & 288 & 312
\end{array}
$$
\end{table}
\begin{table}[htb]
\caption{$N=12m$-point spherical $t$-designs formed from the union of $m$
full orbits under group $[3,3]^+$.
The entries give a list of $m$ orbit representatives $A,B,C$.}
$$
\begin{array}{lrr}
\multicolumn{3}{l}{N=36,~~~t=8} \\
0.74051521 & 0.24352778 & 0.62636367 \\
0.80542549 & 0.30620001 & -0.50747545 \\
0.95712033 & 0.28624872 & 0.04452356 \\
~ \\
\multicolumn{3}{l}{N=60,~~~ t=10} \\
0.71315107 & 0.03408955 & 0.70018102 \\
0.75382867 & 0.54595191 & -0.36562119 \\
0.78335594 & -0.42686412 & -0.45181910 \\
0.93321004 & 0.12033145 & -0.33858436 \\
0.95799794 & 0.27623022 & 0.07705072 \\
~ \\
\multicolumn{3}{l}{N=72,~~~ t=11} \\
0.66932119 & -0.65648669 & -0.34789994 \\
0.75683290 & 0.38164750 & -0.53061205 \\
0.82190371 & 0.54929373 & -0.15083333 \\
0.85544705 & 0.04115447 & 0.51625251 \\
0.90728126 & 0.36233033 & 0.21344190 \\
0.97885492 & 0.12557302 & -0.16147588 \\
~ \\
\multicolumn{3}{l}{N=96,~~~ t=13} \\
0.69989534 & 0.59974524 & -0.38788163 \\
0.73338128 & -0.54971991 & -0.39994990 \\
0.78556905 & 0.09585688 & -0.61130412 \\
0.82321276 & 0.56450535 & 0.06045217 \\
0.83255539 & -0.25643858 & -0.49100996 \\
0.88122889 & 0.33818291 & -0.33025441 \\
0.96391874 & -0.26382492 & -0.03545521 \\
0.96783463 & -0.01683358 & -0.25102343 \\
~ \\
\multicolumn{3}{l}{N=108,~~~ t=14} \\
0.69160471 & -0.40217576 & 0.59994798 \\
0.71050575 & 0.58202818 & 0.39550573 \\
0.75403890 & 0.65127837 & -0.08521631 \\
0.80598041 & 0.26283378 & 0.53039041 \\
0.86226532 & -0.39729017 & 0.31410038 \\
0.86442500 & -0.05628604 & 0.49960114 \\
0.87315060 & -0.46879380 & -0.13356797 \\
0.96418944 & 0.16093133 & 0.21080756 \\
0.97567128 & -0.17376307 & 0.13368600
\end{array}
$$
\end{table}
\setcounter{table}{2}
\begin{table}[htb]
\caption{(cont.) $N=12m$-point spherical $t$-designs formed from the union of $m$
full orbits under group $[3,3]^+$.
The entries give a list of $m$ orbit representatives $A,B,C$.}
$$
\begin{array}{lrr}
\multicolumn{3}{l}{N=144,~~~ t=16} \\
0.65758346 & 0.61920220 & 0.42915339 \\
0.70203400 & -0.68122298 & 0.20756570 \\
0.70428352 & -0.55221381 & 0.44614418 \\
0.71018481 & -0.16518988 & -0.68436090 \\
0.84130836 & -0.32306467 & 0.43339297 \\
0.84532735 & -0.30622774 & -0.43777418 \\
0.85087242 & 0.52354706 & -0.04375603 \\
0.85473787 & -0.02894596 & 0.51825216 \\
0.87135881 & 0.43350173 & 0.22980441 \\
0.94028712 & -0.28839660 & 0.18079695 \\
0.96296114 & 0.02735042 & 0.26824950 \\
0.98473889 & 0.16325742 & 0.06030199 \\
~\\
\multicolumn{3}{l}{N=240,~~~ t=21} \\
0.66536339 & 0.58086027 & 0.46892741 \\
0.67683321 & -0.48257247 & 0.55589623 \\
0.71800639 & 0.65744688 & -0.22853979 \\
0.72687147 & -0.02748828 & -0.68622319 \\
0.73733200 & -0.62085150 & -0.26624225 \\
0.77263286 & 0.51705945 & -0.36835851 \\
0.77909960 & -0.23760971 & -0.58012537 \\
0.78443181 & 0.28431902 & -0.55120724 \\
0.78559925 & -0.40515695 & -0.46763412 \\
0.81763902 & -0.57522572 & 0.02412057 \\
0.84781923 & 0.06632578 & -0.52612113 \\
0.86317647 & -0.46818182 & -0.18902953 \\
0.89265354 & -0.41253405 & 0.18161861 \\
0.89457952 & -0.27876240 & 0.34931219 \\
0.90354264 & 0.09900269 & 0.41690427 \\
0.90950707 & 0.29209374 & 0.29576703 \\
0.94298382 & 0.33269411 & 0.00980574 \\
0.95866803 & -0.10111361 & 0.26595424 \\
0.97946878 & 0.11341985 & 0.16666388 \\
0.99028895 & 0.12883316 & -0.05224764
\end{array}
$$
\end{table}

An explicit formula for $e_t$ (for $t \ge 6$) can be found from \eqn{eq8},
\eqn{eq9}:
\begin{eqnarray*}
\lefteqn{e_t = \left[ \frac{t}{12} \right]
\left(t- 6 \left[ \frac{t}{12} \right] - 5 \right) +
\left[ \frac{t-6}{12} \right] \left ( t-6-6 \left[ \frac{t-6}{12} \right] -5 \right) } \\
&& + 9 \left[ \frac{t}{12} \right] + 9 \left[ \frac{t-6}{12} \right] +
A_{t~\bmod~12} + A_{(t-6) ~\bmod~12} + 1 ~,
\end{eqnarray*}
where $A_0 , \dd, A_{11}$ are 0, 0, 0, 1, 2, 2, 3, 4, 5, 6, 7, 8.
Therefore,
if these designs continue to exist, we will have a sequence of $t$-designs
with $N= 12m$ points satisfying $N= \frac{t^2}{2} (1+ o(1))$ as $t \to \infty$.
Incidentally, Korevaar and Meyers \cite{KM93} show that there exist $t$-designs with $N= O(t^3)$ points and conjecture that $N= O(t^2)$ should be possible.
Eq.~(2) gives a lower bound of $\frac{t^2}{4} (1+ o(1))$.
\paragraph{Acknowledgements.}
We thank Bruce Reznick for some very helpful comments on an earlier version of this paper.

\end{document}